\documentclass[12pt]{amsart}

\newtheorem{theorem}{Theorem}[section]
\newtheorem{lemma}[theorem]{Lemma}
\newtheorem{proposition}[theorem]{Proposition}
\newtheorem{corollary}[theorem]{Corollary}

\newcommand{\adj}{\mathrm{Adj}}
\newcommand{\et}{\quad\mbox{and}\quad}
\newcommand{\bC}{\mathbb{C}}
\newcommand{\bP}{\mathbb{P}}
\newcommand{\bQ}{\mathbb{Q}}
\newcommand{\bR}{\mathbb{R}}
\newcommand{\bZ}{\mathbb{Z}}
\newcommand{\cC}{{\mathcal{C}}}
\newcommand{\cE}{{\mathcal{E}}}
\newcommand{\cK}{{\mathcal{K}}}
\newcommand{\ggll}{\sim}
\newcommand{\GL}{\mathrm{GL}}
\newcommand{\trace}{\mathrm{trace}}

\newcommand{\ud}{\mathbf{d}}
\newcommand{\ue}{\mathbf{e}}

\newcommand{\uu}{\mathbf{u}}
\newcommand{\uun}{\mathbf{1}}
\newcommand{\ux}{\mathbf{x}}
\newcommand{\uy}{\mathbf{y}}
\newcommand{\uz}{\mathbf{z}}
\newcommand{\disp}{\displaystyle}
\newcommand{\dZ}[1]{\left\{#1\right\}}
\newcommand{\Res}{{\mathrm{\,Res}}}

\begin{document}

\baselineskip=17pt

\title{Diophantine approximation in small degree}
\author{Damien ROY}
\address{
   D\'epartement de Math\'ematiques\\
   Universit\'e d'Ottawa\\
   585 King Edward\\
   Ottawa, Ontario K1N 6N5, Canada}
\email{droy@uottawa.ca}
\thanks{Work partially supported by NSERC and CICMA}

\maketitle


\section{Introduction}
\label{sec-intro}


This paper (partly a survey) deals with the problem of finding
optimal exponents in Diophantine estimates involving one real
number $\xi$.  The prototype of such an estimate is the fact,
known at least since Euler, that, for any given irrational real
number $\xi$, there exist infinitely many rational numbers $p/q$
with
\begin{equation}
 \label{rational-estimate}
\left| \xi-{p\over q} \right| \le q^{-2}.
\end{equation}
Here, the exponent of $q$ in the upper bound is optimal because,
when $\xi$ has bounded partial quotients, there is also a constant
$c>0$ such that $|\xi-p/q|\ge cq^{-2}$ for all rational numbers
$p/q$ (see Chapter I of \cite{Sc}).

Define the {\it height} $H(P)$ of a polynomial $P\in\bR[T]$ as the
largest absolute value of its coefficients, and the {\it height}
$H(\alpha)$ of an algebraic number $\alpha$ as the height of its
irreducible polynomial over $\bZ$.  Then the above estimate may be
generalized in the following two ways related respectively with
Mahler's and Koksma's classifications of numbers.

Consider a real number $\xi$ which, for a fixed integer $n\ge 1$,
is not algebraic over $\bQ$ of degree $\le n$.  On one hand, an
application of Dirichlet's box principle shows that there exist
infinitely many non-zero polynomials $P\in \bZ[T]$ such that
\begin{equation}
 \label{polvalue}
 |P(\xi)| \ll H(P)^{-n}
\end{equation}
where, as the sequel, the implied constant depends only on $n$ and
$\xi$.   On the other hand, Wirsing showed in \cite{Wi} that there
exist infinitely many algebraic numbers $\alpha$ of degree at most
$n$ with
\begin{equation}
 \label{Wirsing-estimate}
 |\xi-\alpha| \ll H(\alpha)^{-(n+3)/2}.
\end{equation}

For $n=1$, both estimates are equivalent to
(\ref{rational-estimate}), up to the values of the implied
constants. For general $n$, Spind\v{z}uk proved in \cite{Sp} that
the exponent of $H(P)$ in the upper bound (\ref{polvalue}) is
optimal by showing that, for $\xi$ outside of a set of Lebesgue
measure zero and for each $\epsilon>0$, there are only finitely
many non-zero integer polynomials of degree at most $n$ with
$|P(\xi)|\le H(P)^{-n-\epsilon}$.  However, in the second estimate
(\ref{Wirsing-estimate}), it was conjectured by Wirsing (p.\ 69 of
\cite{Wi}) and Schmidt (p.\ 258 of \cite{Sc}) that the optimal
exponent for $H(\alpha)^{-1}$ is $n+1$ instead of $(n+3)/2$. Aside
from the case $n=1$, this is known to be true only for $n=2$
thanks to work of Davenport and Schmidt \cite{DSa}. Despite of
several refinements by Bernik, Tishchenko and Wirsing, the optimal
exponent remains unknown for any $n\ge 3$ (see Chapter 3 of
\cite{Bu} for more details and references).

In 1969, Davenport and Schmidt \cite{DSb} devised a new method
based on geometry of numbers to study the second type of estimate.
Its flexibility is such that it allowed them to treat
approximation by algebraic integers.  Assuming, for a fixed $n\ge
2$, that $\xi$ is not algebraic over $\bQ$ of degree $\le n-1$,
they proved that there exists infinitely many algebraic integers
$\alpha$ of degree at most $n$ with
\begin{equation}
 \label{DS-estimate}
 |\xi-\alpha| \ll H(\alpha)^{-[(n+1)/2]}.
\end{equation}
They showed that, for $n=2$, the optimal exponent of approximation
is $2$ in agreement with the natural conjecture that it should, in
general, be $n$.  They also provided sharper estimates for
$n=3,4$.

Their approach which we will describe in the next section conveys
them to establish first another Diophantine approximation result.
In the case $n=3$, it is the following statement (Theorem 1a of
\cite{DSb}) where
\[
\gamma=\frac{1+\sqrt{5}}{2}
\]
denotes the golden number.

\begin{theorem} [Davenport-Schmidt]
\label{DSconvex}
Suppose that $\xi\in\bR$ is neither rational nor quadratic over
$\bQ$. Then there are arbitrarily large values of $X$ such that
the inequalities
\begin{equation}
 \label{convex}
 |x_0| \le X,
 \quad
 |x_0\xi-x_1| \le cX^{-1/\gamma},
 \quad
 |x_0\xi^2-x_2| \le cX^{-1/\gamma}
\end{equation}
where $c$ is a suitable positive number depending on $\xi$, have
no solution in integers $x_0$, $x_1$, $x_2$, not all $0$.
\end{theorem}

Note that, an application of Dirichlet's box principle shows that,
for any $X\ge 1$, there exists a non-zero point $(x_0,x_1,x_2)$ in
$\bZ^3$ with $|x_0|\le X$ and $|x_0\xi^j-x_j|\le [\sqrt{X}]^{-1}$
for $j=1,2$ (see Theorem 1A in Chapter II of \cite{Sc}). Since
$1/\gamma\simeq 0.618 >1/2$, the condition (5) is a far stronger
requirement.  Nevertheless, it is shown in \cite{Ra, Rb} that the
exponent $1/\gamma$ is best possible for this problem.  More
precisely, there are countably many real numbers $\xi$ which are
neither rational nor quadratic over $\bQ$ such that, for a
different choice of $c>0$ (depending on $\xi$), the inequalities
(\ref{convex}) admit a non-zero integer solution for any $X\ge 1$
(Theorem 1.1 of \cite{Rb}).  Because of this property, we shall
call these numbers {\it extremal}.

In Section 3 below, we will sketch a proof of Theorem
\ref{DSconvex} and of a criterion for a real number to be
extremal. This criterion attaches to an extremal real number a
sequence of approximation triples which we will show in Section 4
to be essentially unique and to satisfy a certain recurrence
relation.  This will allow us in Section 5 to derive a
construction of extremal real numbers which generalizes that of
Section 6 of \cite{Rb}.

Theorem~\ref{DSconvex} concerns simultaneous approximations of a
real number and its square by rational numbers with the same
denominator.  Davenport and Schmidt looked more generally at
simultaneous approximations of the first $n-1$ powers
$\xi,\dots,\xi^{n-1}$ of a real number $\xi$ by rational numbers
with the same denominator (Theorem 2a of \cite{DSb}), and their
result was recently improved by M.~Laurent \cite{La}. However, the
optimal exponent for this problem is unknown for $n\ge 4$.

\begin{corollary} [Davenport-Schmidt]
 \label{cor-alg3}
Let $\xi$ be as above.  There exist infinitely many algebraic
integers $\alpha$ of degree at most $3$ with
\[
|\xi-\alpha| \ll H(\alpha)^{-\gamma-1}.
\]
\end{corollary}

It is shown in \cite{Rc} that here also the exponent
$\gamma+1\simeq 2.618$ is best possible, against the natural
conjecture that the optimal exponent would be $3$.  More
precisely, there are real numbers for which the above corollary is
optimal up to the value of the implied constant (Theorem 1.1 of
\cite{Rc}).  Such numbers have to be extremal and it would be
interesting to know if this property extends to all extremal real
numbers.

Considering the irreducible polynomials in $\bZ[T]$ of the
approximations $\alpha$ provided by the above corollary, we
readily deduce:

\begin{corollary}
 \label{cor-pol3}
Let $\xi$ be as above.  There exist infinitely many monic
polynomials $P\in\bZ[T]$ of degree at most $3$ with
\[
|P(\xi)| \ll H(P)^{-\gamma}.
\]
\end{corollary}

We will prove in Section 6 that the exponent $\gamma$ in this
statement is also best possible. Note that an argument of Bugeaud
and Teuli\'e \cite{BT, Te} shows more precisely that the
inequality (\ref{DS-estimate}) has infinitely many solutions in
algebraic integers $\alpha$ of degree exactly $n$ under the same
assumption that $\xi$ is not algebraic of degree $\le n-1$.
Therefore, one may require that the algebraic integers of
Corollary \ref{cor-alg3} have degree $3$ and that the polynomials
of Corollary~\ref{cor-pol3} also have degree $3$.

Denote by $\bZ[T]_{\le n}$ the group of integer polynomials of
degree at most $n$.  In another direction, we have the following
Gel'fond type criterion in degree two \cite{AR} which is, in a
sense, dual to Theorem \ref{DSconvex}.

\begin{theorem} [Arbour-Roy]
 \label{ARgelfond}
Let $\xi \in \bC$. Assume that for any sufficiently large positive
real number $X$ there exists a non-zero polynomial $P \in
\bZ[T]_{\le 2}$ of height at most $X$ such that
\begin{equation}
 \label{gelfond-estimate}
        |P(\xi)| \leq \frac{1}{4}X^{-\gamma-1}.
\end{equation}
Then $\xi$ is algebraic over $\bQ$ of degree at most $2$.
\end{theorem}

Theorem 1.2 of \cite{Rb} shows that the exponent of $X$ in
(\ref{gelfond-estimate}) is best possible. Analog statements
involving polynomials of degree at most $n$ for a fixed integer
$n\ge 3$ are well-known but the corresponding optimal exponent is
not known (see Theorem 1 of \cite{Br} for a general setting, and
Theorem 2b of \cite{DSb} for a sharper estimate in the present
context).

Finally, we will show in Section 7 that the extremal real numbers
(associated with Theorem~\ref{DSconvex}) are also characterized as
those real numbers $\xi$ for which Theorem~\ref{ARgelfond} is
optimal up to the value multiplicative constant.

\begin{theorem}
 \label{thm-equiv}
Let $\xi$ be an real number which is not rational nor quadratic
over $\bQ$.  The following conditions are equivalent:
\begin{itemize}
\item[(a)]
 there exists a constant $c>0$ such that, for any real number
 $X\ge 1$, there is a non-zero point $\ux=(x_0,x_1,x_2)\in\bZ^3$
 satisfying the condition (\ref{convex});
\item[(b)]
 there exists a constant $c>0$ such that, for any real number
 $X\ge 1$, there is a non-zero polynomial $P\in\bZ[T]_{\le 2}$
 of height at most $X$ satisfying $|P(\xi)|\le cX^{-\gamma-1}$.
\end{itemize}
\end{theorem}

It would be interesting to know if a similar property holds in
higher degree.  Note that any real number satisfying one of the
above conditions (a) or (b) is transcendental over $\bQ$ by virtue
of Schmidt's subspace theorem.

The present work originates from a colloquium talk given at the
University of Ottawa in November 2002.  The author thanks the
editors for suggesting to include it in the proceedings of CNTA-7.
He also thanks Yann Bugeaud for pointing out a number of
stimulating questions in relation to the present topic, in
particular concerning the status of the exponent in Corollary
\ref{cor-pol3}.

\section{The method of Davenport and Schmidt}

Let $n$ be a fixed positive integer and let $\xi$ be a fixed real
number which is not algebraic over $\bQ$ of degree at most $n$.
The natural approach to construct algebraic approximations of
$\xi$ of degree at most $n$ is to produce non-zero polynomials of
$\bZ[T]_{\le n}$ with ``small'' value and ``large'' derivative at
$\xi$. One then concludes using the fact that any non-zero
polynomial $P\in\bR[T]$ of degree at most $n$ with $P'(\xi)\neq 0$
has at least one root $\alpha$ with
\begin{equation}
  \label{distance-xi-alpha}
  |\xi-\alpha| \le n\frac{|P(\xi)|}{|P'(\xi)|}.
\end{equation}

Define a {\it convex body} of $\bR^n$ to be a compact, convex,
neighborhood $\cC$ of $0$ which is symmetric with respect to $0$
(i.e.\ $\cC=-\cC$).  According to a well-known result of
Minkowski, if such a convex body $\cC$ has volume at least $2^n$,
then it contains a non-zero point of $\bZ^n$.  Applying this
result to the convex body of $\bR^{n+1}$ defined, for a real
number $X\ge 1$, by
\[
 \left\{
 \begin{array}{l}
 |x_0+x_1\xi+\dots+x_n\xi^n| \le X^{-n} \\
 |x_1| \le X \\
 \dots \\
 |x_n| \le X \\
 \end{array}
 \right.
\]
and noting that its volume is $2^{n+1}$, we deduce that there
exists a non-zero polynomial $P\in\bZ[T]$ with
$$
\deg(P)\le n,
  \quad
H(P) \le (1+|\xi|+\dots+|\xi|^n) X
  \et
|P(\xi)| \le X^{-n}.
$$
The difficulty is to control the derivative of $P$ at $\xi$.  The
best one can hope is, for arbitrarily large values of $X$, to have
$|P'(\xi)|\gg X$. Then, using (\ref{distance-xi-alpha}), one finds
that there is at least one root $\alpha$ of $P$ with
$$
\deg(\alpha)\le n,
  \quad
H(\alpha)\ll X
  \et
|\xi-\alpha| \ll X^{-n-1} \ll H(\alpha)^{-n-1}.
$$
This motivates the conjecture mentioned in the introduction.  In
general, one has recourse to resultants to establish lower bounds
on $|P'(\xi)|$.

The approach of Davenport and Schmidt in \cite{DSb} is different
as they require a set $\{P_1,\dots,P_n\}$ of $n$ linearly
independent polynomials of $\bZ[T]_{\le n-1}$, all having small
absolute value at $\xi$. Then, taking any monic polynomial
$Q\in\bR[T]$ of degree $n$ with $Q(\xi)=0$ and writing it as a
linear combination
\[
Q = T^n + \theta_1P_1+\cdots+\theta_nP_n
\]
with $\theta_1,\dots,\theta_n\in \bR$, one sees that the
polynomial $P\in\bZ[T]$ given by
\[
P= T^n +[\theta_1]P_1+\dots+[\theta_n]P_n
\]
is monic of degree $n$, has height $H(P) \le H(Q)+\sum_{i=1}^n
H(P_i)$, and satisfies
\[
 |P(\xi)| \le \sum_{i=1}^n |P_i(\xi)|
 \quad\hbox{as well as}\quad
 |P'(\xi)| \ge |Q'(\xi)| - \sum_{i=1}^n |P'_i(\xi)|.
\]
So, if $X$ denotes an upper bound for $\sum_{i=1}^n H(P_i)$, an
appropriate choice of $Q$ (with $H(Q)\ll X$ and $|Q'(\xi)|\gg X$)
produces a polynomial $P$ with
\[
 H(P)\ll X,
 \quad
 |P(\xi)|\le \sum_{i=1}^n |P_i(\xi)|
 \et
 |P'(\xi)| \ge X.
\]
The roots of such a polynomial are algebraic integers of degree
$\le n$ and height $\ll X$ and, by (\ref{distance-xi-alpha}), at
least one of them, say $\alpha$, satisfies
\[
|\xi-\alpha| \le \frac{n}{X}\sum_{i=1}^n |P_i(\xi)|.
\]

To construct appropriate sets of polynomials $\{P_1,\dots,P_n\}$,
Davenport and Schmidt apply a result of duality of Mahler
\cite{Ma}. To state this result or rather a consequence of it, let
$\cC$ be a convex body of $\bR^n$ and let $\cC^*$ denote the set
of points $(y_1,\dots,y_n) \in \bR^n$ satisfying
\[
|x_1y_1+\cdots+x_ny_n| \le 1
\]
for all $(x_1,\dots,x_n)\in\cC$.  Then, $\cC^*$ is again a convex
body of $\bR^n$, called the {\it dual} (or {\it polar}) convex
body to $\cC$ (the bi-dual $\cC^{**}$ being $\cC$ itself) and, if
$\cC$ contains no non-zero integral point, then $n!\cC^*$ contains
$n$ linearly independent points of $\bZ^n$.

For example, given real numbers $X,Y \ge 1$, the convex bodies of
$\bR^n$ defined by
\[
\cC\ : \left\{
\begin{array}{l}
 |x_0|\le X \\
 |x_0\xi-x_1|\le Y^{-1} \\
 \dots\\
 |x_0\xi^{n-1}-x_{n-1}|\le Y^{-1}\\
\end{array}
\right.
 \et
\cK\ : \left\{
\begin{array}{l}
 |y_0+y_1\xi+\dots+y_{n-1}\xi^{n-1}| \le X^{-1} \\
 |y_1| \le Y \\
 \dots \\
 |y_{n-1}| \le Y \\
\end{array}
\right.
\]
are essentially dual to each other in the sense that
\[
n^{-1}\cK\subseteq \cC^* \subseteq \cK.
\]
So, if $\cC$ contains no non-zero integral point, then there are
$n$ linearly independent polynomials of $\bZ[T]_{\le n-1}$ of
height $\ll Y$ whose absolute values at $\xi$ are $\ll X^{-1}$ and
therefore there exists an algebraic integer of degree $\le n$ and
height $\ll Y$ whose distance to $\xi$ is $\ll (XY)^{-1}$.

In the case $n=3$, this explains why Theorem~\ref{DSconvex}
implies the existence of infinitely many algebraic integers
$\alpha$ of degree $\le 3$ with $|\xi-\alpha| \ll
H(\alpha)^{-\gamma-1}$ as stated in Corollary~\ref{cor-alg3}. In
general, the fact that (\ref{DS-estimate}) has infinitely many
solutions in algebraic integers of degree $\le n$ follows from a
similar statement (Theorem 2a of \cite{DSb}) showing that for $X =
cY^\nu$ with $\nu=[(n-1)/2]$ and an appropriate constant $c>0$,
there are arbitrarily large values of $Y$ for which the convex
body $\cC$ contains no non-zero integral point.

A modification of the method produces approximation by algebraic
numbers or algebraic integers of degree $n$ or even by algebraic
units of degree $n$ (see \cite{BT, Te}). A more general choice of
convex bodies involving several derivatives still produces
simultaneous approximations of a real number by conjugate
algebraic integers \cite{RW}.

\section{Extremal real numbers}

In this section, we present a sketch of proof of Theorem
\ref{DSconvex} and establish some properties of the corresponding
``extremal'' real numbers.

Let $\xi$ be a fixed real number.  For each point
$\ux=(x_0,x_1,x_2)\in\bZ^3$, we define
\[
\|\ux\|=\max\{|x_0|,|x_1|,|x_2|\}
  \et
L(\ux)=L_\xi(\ux)=\max\{|x_1-\xi x_0|, |x_2-\xi^2 x_0|\}.
\]
Identifying any such point with the corresponding symmetric matrix
\begin{equation}
 \label{matrix}
 \ux=\begin{pmatrix} x_0 &x_1\\ x_1 &x_2\\ \end{pmatrix},
\end{equation}
we define
\[
\det(\ux)
  = \left|
    \begin{matrix}
          x_0 &x_1\\ x_1 &x_2\\
    \end{matrix}
    \right|
  = x_0x_2-x_1^2.
\]
Then, using the multilinearity of the determinant, one finds
\begin{equation}
\label{det1}
 |\det(\ux)|
   = \left\|
     \begin{matrix}
     x_0 &x_1-\xi x_0\\ x_1 &x_2-\xi x_1\\
     \end{matrix}
     \right\|
  \ll \|\ux\| L(\ux).
\end{equation}
Similarly, if $\det(\ux,\uy,\uz)$ denotes the determinant of the
$3\times 3$ matrix whose rows are points $\ux, \uy, \uz \in
\bZ^3$, one finds,
\begin{equation}
 \label{det3}
 \begin{array}{rl}
  |\det(\ux,\uy,\uz)|
   &= \left\|
      \begin{matrix}
      x_0 &x_1-\xi x_0 &x_2-\xi^2 x_0\\
      y_0 &y_1-\xi y_0 &y_2-\xi^2 y_0\\
      z_0 &z_1-\xi z_0 &z_2-\xi^2 z_0\\
      \end{matrix}
      \right\| \\ \\
   &\ll \|\ux\| L(\uy) L(\uz)
      + \|\uy\| L(\ux) L(\uz)
      + \|\uz\| L(\ux) L(\uy). \\
 \end{array}
\end{equation}

\medskip
We are now ready to present a sketch of proof of Theorem
\ref{DSconvex}.  To this end, assume that $\xi$ is neither
rational nor quadratic over $\bQ$ and that there exists a positive
real number $c$ such that the inequalities (\ref{convex}) have a
non-zero solution $\ux = (x_0,x_1,x_2) \in \bZ^3$ for any
sufficiently large real number $X$. We need to show that $c$ is
bounded below by some positive constant depending only on $\xi$.

\medskip
First note that there is a sequence of points $(\ux_i)_{i\ge 1}$
in $\bZ^3$ with the following three properties:
\begin{itemize}
\item $1\le \|\ux_1\| < \|\ux_2\| < \|\ux_3\| < \dots$
\item $L(\ux_1) > L(\ux_2) > L(\ux_3) > \dots$
\item if $\uy\in\bZ^3$ has $1\le \|\uy\| < \|\ux_{i+1}\|$, then
  $L(\uy)\ge L(\ux_i)$.
\end{itemize}
Although this differs slightly from the construction in \S3 of
\cite{DSb}, we say that $(\ux_i)_{i\ge 1}$ is a sequence of {\it
minimal points} for $\xi$.  In such a sequence, any point $\ux_i$
with $i\ge 2$ is {\it primitive}, i.e.\ has relatively prime
coordinates. Moreover any two consecutive points $\ux_i$ and
$\ux_{i+1}$ with $i\ge 2$ are linearly independent.  Define
\[
 X_i=\|\ux_i\| \et L_i=L(\ux_i)
\]
for each $i\ge 1$.  Then, for any sufficiently large $i$, the
hypotheses imply the existence of a non-zero point $\ux \in \bZ^3$
with $\|\ux\|<X_{i+1}$ and $L(\ux)\ll c X_{i+1}^{-1/\gamma}$. This
gives
\begin{equation}
 \label{borneL}
 L_i \ll c X_{i+1}^{-1/\gamma}.
\end{equation}
Davenport and Schmidt show that, for any sufficiently large $i$,
we have $\det(\ux_i)\neq 0$ (Lemma 2 of \cite{DSb}).  Since
$\det(\ux_i)$ is an integer, its absolute value is then bounded
below by $1$ and, using the estimate (\ref{det1}) combined with
(\ref{borneL}), we find
\begin{equation}
  \label{borne1}
  1 \le |\det(\ux_i)|
    \ll X_iL_i
    \ll cX_i X_{i+1}^{-1/\gamma}.
\end{equation}
They also show that, for infinitely many $i$, the points
$\ux_{i-1}$, $\ux_{i}$ and $\ux_{i+1}$ are linearly independent
(Lemma 5 of \cite{DSb}). For these $i$, the same argument based on
(\ref{det3}) and (\ref{borneL}) then gives
\begin{equation}
  \label{borne3}
  1 \le |\det(\ux_{i-1}, \ux_i,\ux_{i+1})|
    \ll X_{i+1}L_iL_{i-1}
    \ll c^2 X_{i+1}^{1-1/\gamma}X_i^{-1/\gamma}.
\end{equation}
The required lower bound on $c$ then follows by comparing
(\ref{borne1}) and (\ref{borne3}), upon noting that
$1-1/\gamma=1/\gamma^2$.

\medskip
The above considerations apply in particular to any extremal real
number $\xi$.  Combining (\ref{borne1}) and (\ref{borne3}) then
shows that, for all indices $i$ such that $\ux_{i-1}$, $\ux_{i}$
and $\ux_{i+1}$ are linearly independent, we have
\[
 \|\ux_{i+1}\| \ggll \|\ux_i\|^\gamma,
 \quad
 L_\xi(\ux_i) \ggll \|\ux_i\|^{-1},
 \quad
 |\det(\ux_i)|\ggll 1,
 \quad
 |\det(\ux_{i-1},\ux_i,\ux_{i+1})|\ggll 1
\]
writing $a\ggll b$ to mean $a\ll b$ and $a\gg b$.  A further
analysis shows that, by going to a subsequence $(\uy_i)_{i\ge 1}$
of $(\ux_i)_{i\ge 1}$, one may assume these estimates to hold for
all $i\ge 1$.  More precisely, we have the following equivalence
(Theorem 5.1 of \cite{Rb}):

\begin{theorem}
   \label{thm-extremal}
A real number $\xi$ is extremal if and only if there exists a
constant $c\ge 1$ and an unbounded sequence of non-zero primitive
points $(\uy_i)_{i\ge 1}$ of $\bZ^3$ satisfying, for all $i\ge 1$,
\begin{enumerate}
\item[]
  $c^{-1} \|\uy_{i}\|^\gamma
    \le \|\uy_{i+1}\|
    \le c \|\uy_{i}\|^\gamma$,
\item[]
  $c^{-1} \|\uy_i\|^{-1} \le L_\xi(\uy_i) \le c \|\uy_i\|^{-1}$,
\item[]
  $1\le |\det(\uy_i)| \le c$,
\item[]
  $1\le |\det(\uy_i,\uy_{i+1},\uy_{i+2})| \le c$.
\end{enumerate}
\end{theorem}

In the next section, we show that the sequence $(\uy_i)_{i\ge 1}$
is essentially uniquely determined by $\xi$.

\section{The sequence of approximation triples
\label{apptriples}}

In this section, we fix an extremal real number $\xi\in\bR$ and a
sequence of approximation triples $(\uy_i)_{i\ge 1}$ as in
Theorem~\ref{thm-extremal}.  We first prove:

\begin{proposition}
There exists a constant $c_3>0$ such that any non-zero primitive
point $\uy\in\bZ^3$ with
\begin{equation}
  \label{conditiony}
  L_\xi(\uy) \le c_3 \|\uy\|^{-1/\gamma}
\end{equation}
is of the form $\uy=\pm \uy_i$ for some index $i\ge 1$.
\end{proposition}

\begin{proof}
Fix a constant $c_3$ with $0<c_3 \le \|\uy_1\|^{-1}$.  Then, for
any non-zero point $\uy\in\bZ^3$, there exists an index $i\ge 1$
such that
\[
c_3 \|\uy_i\| \le \|\uy\| < c_3 \|\uy_{i+1}\|.
\]
If (\ref{conditiony}) holds, we then find
\[
\begin{array}{rclcl}
|\det(\uy,\uy_i,\uy_{i+1})|
  &\ll &\|\uy\| \|\uy_i\|^{-1} \|\uy_{i+1}\|^{-1}
     + \|\uy_{i+1}\| \|\uy_i\|^{-1} L_\xi(\uy)
  &\ll &c_3 + c_3^{1/\gamma^2}, \\
|\det(\uy,\uy_i,\uy_{i-1})|
  &\ll &\|\uy\| \|\uy_i\|^{-1} \|\uy_{i-1}\|^{-1}
     + \|\uy_i\| \|\uy_{i-1}\|^{-1} L_\xi(\uy)
  &\ll &c_3 +c_3^{1/\gamma^2}. \\
\end{array}
\]
So, provided that $c_3$ is sufficiently small, these determinants
vanish and, since $\uy_{i-1}$, $\uy_i$ and $\uy_{i+1}$ are
linearly independent with $\uy_i$ primitive, we conclude that
$\uy=\pm\uy_i$.
\end{proof}

Since, for all sufficiently large values of $i$, the point
$\uy=\uy_i$ satisfies the condition (\ref{conditiony}), we deduce
from this proposition that the sequence $(\uy_i)_{i\ge 1}$ is
uniquely determined by $\xi$ up to its first terms and up to
multiplication of its terms by $\pm 1$.

\medskip
In proving the inequalities (\ref{det1}) and (\ref{det3}), we used
the multi-linearity of the determinant.  Equivalently, we could
have looked at the Taylor series expansion of $\det(\ux)$ and
$\det(\ux,\uy,\uz)$ at the points $\ux=(x_0,x_0\xi,x_0\xi^2)$,
\dots, $\uz=(z_0,z_0\xi,z_0\xi^2)$.  Proposition~\ref{searchpol}
of the appendix generalizes this idea and, through a computer
search, provided the following relations.

\begin{proposition}
  \label{prop-rel}
For any sufficiently large index $i$, we have
\begin{equation}
 \label{relations}
 \det(\uy_i,\uy_{i+1},[\uy_{i+3},\uy_{i+3},\uy_{i+4}]) = 0
 \et
 \det(\uy_{i+1},\uy_{i+2},[\uy_{i+3},\uy_{i+3},\uy_{i+4}]) = 0,
\end{equation}
where, upon identifying points $\ux,\uz\in\bZ^3$ with the
corresponding symmetric matrices as in (\ref{matrix}) and upon
denoting by $\adj(\uz)$ the adjoint matrix of $\uz$, we define
\begin{equation}
 \label{bracket}
 [\ux,\ux,\uz] = \ux\, \adj(\uz)\, \ux.
\end{equation}
\end{proposition}

A direct proof of these relations can be found in \cite{Rb}, as
part of the proof of Corollary 5.2 of \cite{Rb}.  It uses the
estimates of Lemma 3.1 of \cite{Rb} to show that the above
determinants (\ref{relations}) have absolute values tending to
zero as $i$ tends to infinity.  As in Proposition 2.3 of
\cite{Rc}, we deduce:

\begin{corollary}
There exists a $2\times 2$ matrix $M$ with relatively prime
integer coefficients and an index $i_0$ such that the symmetric
matrix corresponding to $\uy_{i+2}$ is a rational multiple of
$\uy_{i+1}M\uy_i$ when $i\ge i_0$ is odd, and a rational multiple
of $\uy_{i+1}{^t}M\uy_i$ when $i\ge i_0$ is even.  Such a matrix
$M$ is non-singular, non-symmetric and non-skew-symmetric.
\end{corollary}

\begin{proof}
Choose $i_0\ge 2$ so that the relations (\ref{relations}) hold for
$i\ge i_0-1$. Since $\uy_i$, $\uy_{i+1}$ and $\uy_{i+2}$ are
linearly independent, these relations imply that
$[\uy_{i+3},\uy_{i+3},\uy_{i+4}]$ is a rational multiple of
$\uy_{i+1}$ and thus, by definition of the latter symbol, since
all these matrices are invertible, that $\uy_{i+4}$ is a rational
multiple of $\uy_{i+3} \uy_{i+1}^{-1} \uy_{i+3}$ for $i\ge i_0-1$.
Therefore, assuming that $\uy_{i+3}$ is a rational multiple of
$\uy_{i+2}S\uy_{i+1}$ for some integer matrix $S$ and some index
$i\ge i_0-1$, we find that $\uy_{i+4}$ is a rational multiple of
$\uy_{i+2}S\uy_{i+3}$, and so, by taking transpose, that
$\uy_{i+4}$ is a rational multiple of $\uy_{i+3}{^t}S\uy_{i+2}$.
The first assertion of the corollary then follows by induction on
$i$, upon choosing $M$ so that it holds for $i=i_0$.  The matrix
$M$ is clearly non-singular.  It is not symmetric since a simple
computation based for example on the formulas (2.1) and (2.2) of
\cite{Rb} gives
\begin{equation}
  \label{traceMJ}
  \det(\uy_i,\uy_{i+1},\uy_{i+1}M\uy_i)
  =
  \det(\uy_i)\det(\uy_{i+1})\trace(MJ)
  \quad
  \hbox{where}
  \quad
  J = \begin{pmatrix} 0&1\\ -1&0\\ \end{pmatrix},
\end{equation}
while, for odd $i\ge i_0$, the above determinant is non-zero.
Finally, it is not skew-symmetric, otherwise in the above notation
we would have $M=\pm J$ which, for $i\ge i_0$, would imply
proportionality relations
\[
 \uy_{i+3}
 \propto \uy_{i+2} J \uy_{i+1}
 \propto \uy_i J \uy_{i+1} J \uy_{i+1}
 \propto \uy_i
\]
and thus $\uy_{i+3}=\pm\uy_i$, against the fact that the norms of
the $\uy_i$'s are unbounded.
\end{proof}

As the sequence $(\uy_i)_{i\ge 1}$ is uniquely determined by $\xi$
up to its first terms and up to multiplication of its terms by
$\pm 1$, we deduce that the matrix $M$ of the corollary is
uniquely determined by $\xi$ up to multiplication by $\pm 1$ and
up to transposition. We say that $M$ is the matrix {\it
associated} with $\xi$ and, for a given $M$, we denote by $\cE(M)$
the set of extremal real numbers with associated matrix $M$. In
the next section, we present a criterion for showing that $\cE(M)$
is not empty.

Remark that, in the notation of the corollary, the sequence of
matrices $(M_i)_{i\ge 1}$ given by $M_i=\uy_i M$ for $i$ even and
by $M_i=\uy_i{^t}M$ for $i$ odd satisfies
\[
 M_{i+2}=M_{i+1}M_i
\]
for all $i\ge i_0$.  It can therefore be viewed as a Fibonacci
sequence of matrices in $\GL_2(\bQ)$.

\section{Construction of extremal real numbers}

\begin{proposition}
 \label{criterion}
 Let $M$ be a non-singular, non-symmetric
$2\times 2$ matrix with relatively prime integer coefficients.
Assume that there exist non-singular $2\times 2$ symmetric
matrices $\uy_1,\uy_2,\uy_3$ with relatively prime integer
coefficients such that $\uy_3$ is a rational multiple of
$\uy_2M\uy_1$. Extend the definition of $\uy_i$ coherently for
$i\ge 4$ by asking that $\uy_i$ has relatively prime integer
coefficients and that it is a rational multiple of
$\uy_{i-1}M\uy_{i-2}$ for odd $i\ge 3$ and a rational multiple of
$\uy_{i-1}{^t}M\uy_{i-2}$ for even $i\ge 4$. Assume further that
the $\uy_i$'s are unbounded and that there exist positive
constants $c_4$ and $c_5$ such that
\[
 |\det(\uy_i)| \le c_4
 \et
 \|\uy_{i+2}\| \ge c_5 \|\uy_{i+1}\| \|\uy_i\|
\]
for all $i\ge 1$.  Then $(\uy_i)_{i\ge 1}$ is a sequence of
approximation triples associated with an extremal real number
$\xi\in\cE(M)$.
\end{proposition}

\begin{proof}
For each $i\ge 1$, we have
\begin{equation}
 \label{recurrence}
 \uy_{i+2} = \rho_i \uy_{i+1}S\uy_i
\end{equation}
for some appropriate choice of $S=M$ or ${^t}M$ and some non-zero
rational number $\rho_i$ with $|\rho_i|\le 1$.  Assuming that
$\uy_i$ and $\uy_{i+1}$ are non-singular, this shows that
$\uy_{i+2}$ is also non-singular. Thus, by induction, all
$\uy_i$'s are non-singular, and so
\[
 1 \le |\det(\uy_i)| \le c_4
\]
for $i\ge 1$.  We also deduce that
\[
 \uy_{i+3}
 = \rho_{i+1} \uy_{i+2} {^t}S \uy_{i+1}
 = \rho_i\rho_{i+1} \uy_{i+1} S \uy_i {^t}S \uy_{i+1}
\]
which, starting from the fact that $\uy_1$, $\uy_2$ and $\uy_3$
are symmetric, implies by induction that all $\uy_i$'s are
symmetric.  The formula (\ref{traceMJ}) with $M$ replaced by $S$
also gives
\[
 \det(\uy_i,\uy_{i+1},\uy_{i+2})
 = \rho_i\det(\uy_i)\det(\uy_{i+1})\trace(SJ).
\]
Since $M$ is non-symmetric, we have $\trace(SJ)=\pm\trace(MJ)\neq
0$ and therefore, the above determinant being an integer, it
satisfies
\[
 1
 \le
 |\det(\uy_i,\uy_{i+1},\uy_{i+2})|
 \le
 c_4^2|\trace(MJ)|.
\]
The relation (\ref{recurrence}) also implies that
\[
 \|\uy_{i+2}\| \le c_6 \|\uy_{i+1}\| \|\uy_i\|
\]
where $c_6$ denotes the sum of the absolute values of the
coefficients of $M$.  Defining $q_i = \|\uy_{i+1}\|
\|\uy_i\|^{-\gamma}$ and using the similar lower bound for
$\|\uy_{i+2}\|$ from the hypotheses of the theorem, we deduce that
\[
 c_5 q_i^{-1/\gamma} \le q_{i+1} \le c_6 q_i^{-1/\gamma},
\]
which, by induction, implies $c_7^{-1} \le q_i\le c_7$ with
$c_7=\max\{q_1,q_1^{-1},c_5^{-\gamma^2},c_6^{\gamma^2}\}$ and so
\[
 c_7^{-1} \|\uy_i\|^\gamma
 \le
 \|\uy_{i+1}\|
 \le
 c_7 \|\uy_i\|^\gamma,
\]
for all $i\ge 1$.

Denote by $[\ux]$ the image of a non-zero point $\ux$ of $\bR^3$
in the projective space $\bP^2(\bR)$, and, for any non-zero point
$\uy\in\bR^3$, define
\[
 d([\ux],[\uy])
 = d(\ux,\uy)
 = \frac{\|\ux\wedge\uy\|}{\|\ux\| \|\uy\|}
\]
where $\ux\wedge\uy$ denotes the vector product of $\ux$ and
$\uy$.  This distance function defines the usual topology on
$\bP^2(\bR)$ (see for example Lemma 1.16 of \cite{Ph}) and it is
easily proved to satisfy
\begin{equation}
 \label{distance}
 d(\ux,\uz) \le d(\ux,\uy) +2d(\uy,\uz)
\end{equation}
for any three non-zero points $\ux,\uy,\uz$ of $\bR^3$.  Since
$\bP^2(\bR)$ is compact, the sequence $([\uy_i])_{i\ge 1}$ has an
accumulation point $[\uy]$ for some non-zero $\uy\in\bR^3$. Since
the points $\uy_i$ have bounded determinant and norm tending to
infinity with $i$, we deduce that, by continuity,
\begin{equation}
 \label{dety}
 \det(\uy) = 0.
\end{equation}

In order to estimate the distance between two consecutives points
of this sequence, we note that
\[
 \uy_{i+2}J\uy_{i+1}
 = \rho_i \uy_{i}{^t}S\uy_{i+1}J\uy_{i+1}
 = \rho_i \det(\uy_{i+1}) \uy_{i}{^t}SJ,
\]
thus
\[
 \max_{k,\ell=0,1}
   |y_{i+2,k}y_{i+1,\ell+1}-y_{i+2,k+1}y_{i+1,\ell}|
 \le c_4 c_6 \|\uy_{i}\|,
\]
and so
\begin{equation}
 \label{distance2}
 d(\uy_{i+1},\uy_{i+2})
 \le
 2 c_4 c_6 \frac{\|\uy_{i}\|}{\|\uy_{i+1}\| \|\uy_{i+2}\|}
 \le
 c_8 \|\uy_{i+1}\|^{-2}
\end{equation}
with $c_8=2c_4c_6/c_5$.  Using (\ref{distance}), (\ref{distance2})
and the fact that the norms of the $\uy_i$'s grow faster than any
geometric series, we deduce that, for $k>i\ge 2$, we have
\[
 d(\uy_i,\uy_k)
 \le
 \sum_{j=i}^{k-1} 2^{j-i} d(\uy_j,\uy_{j+1})
 \le
 c_8 \sum_{j=i}^{k-1} 2^{j-i} \|\uy_j\|^{-2}
 \le
 c_9 \|\uy_{i}\|^{-2}
\]
for some constant $c_9>0$.  As $d(\uy_k,\uy)$ can be made
arbitrarily small for a suitable choice of $k>i$, this implies
\[
 d(\uy_i,\uy) \le c_9\|\uy_i\|^{-2},
\]
showing in particular that the sequence $([\uy_i])_{i\ge 1}$
converges to $[\uy]$.

We claim that the point $\uy=(y_0,y_1,y_2)$ has $y_0\neq 0$.
Otherwise, upon denoting by $k$ a fixed index for which $y_k\neq
0$, we would have, for $i\ge 2$,
\[
 |y_ky_{i,0}|
 =
 |y_ky_{i,0}-y_0y_{i,k}|
 \le
 \|\uy\| \|\uy_i\| d(\uy_i,\uy)
 \le
 c_9 \|\uy\| \|\uy_i\|^{-1}.
\]
As this upper bound tends to zero for $i\to\infty$, this would
force the integer $y_{i,0}$ to be zero for all sufficiently large
values of $i$, against the fact that the determinant of any three
consecutive $\uy_i$'s is non-zero.

Since $y_0\neq 0$, we may assume without loss of generality that
$y_0=1$. Writing $\xi=y_1$, we then deduce, by virtue of
(\ref{dety}), that
\[
 \uy = (1,\xi,\xi^2).
\]
and so
\[
 L_\xi(\uy_i)
 \le
 \|\uy\wedge\uy_i\|
 \le
 \|\uy\| \|\uy_i\| d(\uy_i,\uy)
 \le
 c_9 \max\{1,\xi^2\} \|\uy_i\|^{-1}
\]
for any $i\ge 2$.  We also get a lower bound of the same type for
$ L_\xi(\uy_i)$ by combining the estimate (\ref{det1}) with the
lower bound $|\det(\uy_i)| \ge 1$.  Therefore, by
Theorem~\ref{thm-extremal}, the number $\xi$ is extremal, and
$(\uy_i)_{i\ge 1}$ is an associated sequence of approximation
triples.
\end{proof}

To apply the above proposition for a given skew-symmetric matrix
$M$, one has to choose $\uy_1$ and $\uy_2$ so that
\begin{itemize}
\item
 the product $\uy_2M\uy_1$ is symmetric,
\item
 the $\uy_i$'s have bounded non-zero determinants,
\item
 the $\uy_i$'s are unbounded and the ratios $\|\uy_{i+2}\| /
 (\|\uy_{i+1}\| \|\uy_i\|)$ are bounded below by some positive
 constant.
\end{itemize}
The second condition is automatically fulfilled if $M$, $\uy_1$
and $\uy_2$ have determinant $\pm 1$, because then all $\uy_i$'s
have determinant $\pm 1$.  The third condition is also fulfilled
if, for example, the coefficients of $M$ are positive while those
of $\uy_1$ and $\uy_2$ are non-negative.

\medskip
\noindent {\bf Example 1.} If we define
\[
 A = \begin{pmatrix} a &1\\ 1 &0 \end{pmatrix},
 \quad
 B = \begin{pmatrix} b &1\\ 1 &0 \end{pmatrix}
 \et
 M = AB = \begin{pmatrix} ab+1 &a\\ b &1 \end{pmatrix}
\]
for a choice of distinct positive integers $a$ and $b$, then
\[
 \uy_1 = A,
 \quad
 \uy_2 = ABA
 \et
 \uy_3 = \uy_2 M\uy_1 = ABAABA
\]
are symmetric matrices of determinant $\pm 1$ with non-negative
entries while $M$ has positive entries.  The resulting sequence
$(\uy_i)_{i\ge 1}$ thus fulfills all requirements of Proposition
\ref{criterion}.  It can be shown that the corresponding extremal
real number has continued fraction expansion
\[
 \xi = [0,a,b,a,a,b,a,\dots]
\]
given by the Fibonacci word on $\{a,b\}$ (see Theorem 2.2 of
\cite{Ra} or Corollary 6.3 of \cite{Rb}).

\medskip
\noindent {\bf Example 2.} Take
\[
 \uy_1 = \begin{pmatrix} 1 &1\\ 1 &0 \end{pmatrix},
  \quad
 \uy_2 = \begin{pmatrix} a^3+2a &a^3-a^2+2a-1\\
         a^3-a^2+2a-1 &a^3-2a^2+3a-2 \end{pmatrix}
 \et
 M = \begin{pmatrix} a &1\\ -1 &0 \end{pmatrix}
\]
for a fixed positive integer $a$. One readily checks that $\uy_2 M
\uy_1$ is symmetric, that $\det(\uy_1)=\det(\uy_2)=-1$ and that
$\det(M)=1$.  Thus, to ensure that the corresponding sequence
$(\uy_i)_{i\ge 1}$ defines an extremal real number, it remains
only to verify the growth condition on the norms of these points.
To this end, we note that if $\ux, \uy \in \bZ^3$ have coordinates
satisfying $x_0\ge x_1\ge x_2\ge 0$ and $y_0\ge y_1\ge y_2\ge 0$
and if the product $\uz = \uy M\ux$ is symmetric, then we also
have $z_0\ge z_1\ge z_2\ge 0$ and moreover $z_0 \ge (a-1)y_0x_0$.
By recurrence on $i$, using $\uy_{i+2}=\uy_{i+1}M\uy_i$ for odd
$i$ and $\uy_{i+2}=\uy_i M\uy_{i+1}$ for even $i$, we deduce that
$\|\uy_i\|=y_{i,0}$ and that $\|\uy_{i+2}\| \ge (a-1)
\|\uy_{i+1}\| \|\uy_i\|$ for all $i\ge 1$.  So, for $a\ge 2$, the
required growth condition is satisfied (we have $\lim_{i\to\infty}
\|\uy_i\| = \infty$ since $\|\uy_2\|>\|\uy_1\|=1$). In the case
where $a=1$, one finds that the points $\ux=\uy_3$ and $\uy=\uy_2$
satisfy the stronger conditions $x_0\ge 2x_1\ge 4x_2\ge 0$ and
$y_0\ge 2y_1\ge 4y_2\ge 0$ and that, for such points
$\ux,\uy\in\bZ^3$, when the product $\uz = \uy M\ux$ is symmetric,
we also have $z_0\ge 2z_1\ge 4z_2\ge 0$ and $z_0 \ge y_0x_0/2$. By
recurrence on $i$, this gives $\|\uy_i\|=y_{i,0}$ and
$\|\uy_{i+2}\| \ge (1/2) \|\uy_{i+1}\| \|\uy_i\|$ for all $i\ge
2$. In particular, since $\|\uy_2\|=3$ and $\|\uy_3\|=5$, we
deduce that $\|\uy_{i+1}\|>\|\uy_i\|>2$ for $i\ge 1$ and so
$\lim_{i\to\infty} \|\uy_i\| = \infty$.  Thus, in all cases, the
sequence $(\uy_i)_{i\ge 1}$ defines an extremal real number in
$\cE(M)$. This proves the remark at the end of \S3 of \cite{Rc}.

\section{Approximation by cubic algebraic integers}

In order to show that the exponent in Corollary~\ref{cor-pol3} is
best possible, we apply the following criterion where, for a real
number $x$, the symbol $\{x\}$ denotes the distance from $x$ to a
closest integer (compare with Proposition 9.1 of \cite{Rb}).

\begin{lemma}
Let $\xi$ be an extremal real number and let $(\uy_i)_{i\ge 1}$ be
a corresponding sequence of approximation triples.  Assume that
there exists a constant $c_1>0$ such that
\begin{equation}
 \label{dZ}
 \dZ{y_{i,0}\xi^3} \ge c_1
\end{equation}
for any sufficiently large index $i$.  Then there exists a
constant $c_2>0$ such that, for any monic polynomial $P \in
\bZ[T]_{\le 3}$, we have
\[
 |P(\xi)| \ge c_2 H(P)^{-\gamma}.
\]
\end{lemma}

\begin{proof}
Choose an index $i_0\ge 1$ such that (\ref{dZ}) holds for each
$i\ge i_0$.  Multiplying $P$ by a suitable power of $T$ if
necessary, we may assume without loss of generality that $P$ has
degree three. Writing $P(T)=T^3+pT^2+qT+r$, we find, for any $i\ge
i_0$,
\[
\begin{array}{rl}
\dZ{y_{i,0}\xi^3}
  &\le |y_{i,0}P(\xi)| + |p| \dZ{y_{i,0}\xi^2}
       + |q| \dZ{y_{i,0}\xi} \\
  &\le \|\uy_i\| |P(\xi)| + 2 H(P) L_\xi(\uy_i) \\
  &\le \|\uy_i\| |P(\xi)| + c_3 H(P) \|\uy_i\|^{-1} \\
\end{array}
\]
with a constant $c_3>0$ depending only on $\xi$.  Choosing $i$ to
be the smallest integer $i\ge i_0$ for which
$$
 H(P) \le \frac{c_1}{2c_3} \|\uy_i\|,
$$
and using (\ref{dZ}) this implies
$$
 |P(\xi)| \ge \frac{c_1}{2} \|\uy_i\|^{-1}.
$$
The conclusion follows as the above choice of $i$ implies
$\|\uy_i\| \ll H(P)^\gamma$.
\end{proof}

Let $a\ge 1$ be an integer and, for shortness, define $\cE_a=
\cE\begin{pmatrix} a &1\\ -1 &0\end{pmatrix}$.  Example 2 of the
preceding section shows that this set is not empty.  Moreover, it
is proved in \S4 of \cite{Rc} that any extremal real number $\xi$
in $\cE_a$ satisfies the hypotheses of the preceding lemma (more
precisely it is shown that $\lim_{j\to\infty}
\dZ{y_{i+3j,0}\xi^3}$ exists and is a positive real number for
$i=0,1,2$).  So, we deduce:

\begin{theorem}
Let $a\ge 1$ be an integer and let $\xi\in\cE_a$. There exists a
constant $c>0$ such that, for any monic polynomial $P \in
\bZ[T]_{\le 3}$, we have
$$
|P(\xi)| \ge  c H(P)^{-\gamma}.
$$
\end{theorem}

\begin{corollary}
With $\xi$ as above, there exists a constant $c'>0$ such that, for
any algebraic integer $\alpha$ of degree at most $3$ over $\bQ$,
we have
$$
|\xi-\alpha| \ge  c' H(\alpha)^{-\gamma-1}.
$$
\end{corollary}

\section{Proof of Theorem \ref{thm-equiv}}

Theorem 8.1 of \cite{Rb} provides a sequence of polynomials in
$\bZ[T]_{\le 2}$ with small absolute value at a given extremal
real number $\xi$.  These polynomials are constructed by taking
exterior products of consecutive approximation triples of $\xi$.
As a corollary, this result implies Theorem 1.2 of \cite{Rb} which
in turn tells us that (a) implies (b) in Theorem \ref{thm-equiv}.
Our proof that (b) implies (a) will be similar.  We will first
prove that, for any given real number $\xi$ which satisfies this
condition (b), there is a sequence of polynomials in $\bZ[T]_{\le
2}$ with properties parallel to those stated in Theorem
\ref{thm-extremal}. Then, by taking exterior products of
consecutive polynomials in this sequence, we will get points $\ux$
in $\bZ^3$ for which $L_\xi(\ux)$ is small. This will require the
following lemma.

\begin{lemma}
\label{lemma7-1} Let $\xi$ be a real number and let
$P=p_0+p_1T+p_2T^2$ and $Q=q_0+q_1T+q_2T^2$ be polynomials of
$\bZ[T]_{\le 2}$ with
\[
 2H(P)|Q(\xi)| \le H(Q)|P(\xi)|.
\]
Then the point
\[
 \ux
 =
 P\wedge Q
 =
 (p_2q_1-p_1q_2, p_0q_2-p_2q_0, p_1q_0-p_0q_1) \in \bZ^3
\]
satisfies
\[
 (2\max\{1,|\xi|+|\xi^2|\})^{-1} H(Q)|P(\xi)|
   \le L_\xi(\ux)
   \le \frac{3}{2} H(Q)|P(\xi)|.
\]
\end{lemma}

\begin{proof} The upper bound follows immediately from the relations
\[
 x_0\xi-x_1 = p_2Q(\xi)-q_2P(\xi)
 \et
 x_0\xi^2-x_2 = q_1P(\xi)-p_1Q(\xi).
\]
For the lower bound, we simply note that
\[
\begin{array}{rl}
 q_0P(\xi)-p_0Q(\xi)
 &=-\xi(q_1P(\xi)-p_1Q(\xi))-\xi^2(q_2P(\xi)-p_2Q(\xi)) \\
 &=-\xi(x_0\xi^2-x_2)+\xi^2(x_0\xi-x_1)
\end{array}
\]
implies
\[
 \max\{1,|\xi|+|\xi^2|\} L_\xi(\ux)
 \ge \max_{0\le i\le 2} |q_iP(\xi)-p_iQ(\xi)|
 \ge H(Q)|P(\xi)|-H(P)|Q(\xi)|.
\]
\end{proof}

We now prove Theorem \ref{thm-equiv} by adding one more equivalent
condition (compare with Theorem 5.1 of \cite{Rb}):

\begin{theorem}
   \label{theorem7-1}
Let $\xi$ be a real number. The following conditions are
equivalent:
\begin{itemize}
\item[(a)]
 the number $\xi$ is extremal;
\item[(b)]
 the number $\xi$ is neither rational nor quadratic over $\bQ$
 and, for any real number $X\ge 1$, there is a non-zero polynomial
 $P\in\bZ[T]_{\le 2}$ of height at most $X$ satisfying $|P(\xi)|
 \le c_1 X^{-\gamma-1}$ with a constant $c_1=c_1(\xi)$;
\item[(c)]
 There exists a constant $c_2\ge 1$ and an unbounded sequence
 of non-zero polynomials $(Q_k)_{k\ge 1}$ of $\bZ[T]_{\le 2}$ with
 relatively prime coefficients satisfying, for all $k\ge 1$,
\begin{itemize}
\item[]
  $c_2^{-1} H(Q_k)^\gamma
    \le H(Q_{k+1})
    \le c_2 H(Q_{k})^\gamma$,
\item[]
  $c_2^{-1} H(Q_k)^{-\gamma^3}
   \le |Q_k(\xi)|
   \le c_2 H(Q_k)^{\gamma^3}$,
\item[]
  $1\le |\Res(Q_k,Q_{k+1})| \le c_2$,
\item[]
  $1\le |\det(Q_k,Q_{k+1},Q_{k+2})| \le c_2$.
\end{itemize}
\end{itemize}
\end{theorem}

\begin{proof}
As mentioned earlier, Theorem 1.2 of \cite{Rb} shows that (a)
implies (b).

Assume now that (b) is satisfied.  We prove (c) by going back to
the arguments of \cite{AR}.  First of all, we recall that there is
a sequence of ``minimal polynomials'' $(P_i)_{i\ge 1}$ in
$\bZ[T]_{\le 2}$ with the following three properties:
\begin{itemize}
\item $1\le H(P_1) < H(P_2) < H(P_3) < \dots$
\item $|P_1(\xi)| > |P_2(\xi)| > |P_3(\xi)| > \dots$
\item if $P\in\bZ[T]_{\le 2}$ has $1\le H(P) < H(P_{i+1})$,
  then $|P(\xi)|\ge |P_i(\xi)|$
\end{itemize}
(see \S3 of \cite{DSa} or Lemma 5 of \cite{AR}).  For each $i\ge
1$, let $V_i$ denote the sub-$\bQ$-vector space of $\bQ[T]$
generated by $P_i$ and $P_{i+1}$.  Then, for $i\ge 2$, the
polynomials $P_i$ and $P_{i+1}$ form a basis of the group $V_i\cap
\bZ[T]$ of integral polynomials in $V_i$ (see the proof of Lemma 2
of \cite{DSa}).  In particular, $P_i$ has relatively prime
coefficients for $i\ge 2$.  Moreover the condition (b) implies
\begin{equation}
 \label{proof7-1}
 |P_i(\xi)| \le c_1 H(P_{i+1})^{-\gamma-1}
\end{equation}
for any $i\ge 1$.

Let $I$ denote the set of indices $i\ge 2$ for which $P_{i-1}$,
$P_i$ and $P_{i+1}$ are linearly independent.  Lemma 6 of
\cite{AR} shows that $I$ is an infinite set and the arguments in
\S3 of \cite{AR} show that there exists an index $i_1\in I$ such
that $\Res(P_i,P_{i+1})\neq 0$ for all $i\in I$ with $i\ge i_1$.
We define a sub-sequence $(Q_k)_{k\ge 1}$ of $(P_i)_{i\ge 1}$ by
putting $Q_k=P_{i_k}$ where $i_k$ denotes the $k$-th element $i$
of $I$ with $i\ge i_1$.  We claim that this sequence enjoys all
properties stated in (c).

Take $i=i_k$ for some $k\ge 1$.  Using Lemmas 2 and 4 of \cite{AR}
together with (\ref{proof7-1}) , we find
\[
 1
 \le
 |\det(P_{i-1},P_i,P_{i+1})|
 \le
 6 H(P_i)H(P_{i+1})|P_{i-1}(\xi)|
 \le
 6c_1 H(P_i)^{-\gamma} H(P_{i+1})
\]
and
\[
 1
 \le
 |\Res(P_i,P_{i+1})|
 \le
 12H(P_i)H(P_{i+1})^2 |P_i(\xi)|
 \le
 12c_1H(P_i)H(P_{i+1})^{-1/\gamma}.
\]
Comparing these two sets of inequalities, we deduce
\begin{equation}
 \label{proof7-2}
 H(P_{i+1}) \ggll H(P_i)^\gamma,
 \quad
 |P_{i-1}(\xi)| \ggll H(P_i)^{-\gamma-1}
 \et
 |P_i(\xi)| \ggll H(P_{i+1})^{-\gamma-1}.
\end{equation}
In particular, the polynomial $Q_k=P_i$ satisfies
\begin{equation}
 \label{proof7-3}
 |Q_k(\xi)| \ggll H(Q_k)^{-\gamma^3}.
\end{equation}
Moreover, if $i$ is large enough, the pairs $(P,Q)=(P_{i-1},P_i)$
and $(P,Q)=(P_i,P_{i+1})$ both satisfy the hypotheses of Lemma
\ref{lemma7-1} and so we get
\begin{equation}
\label{proof7-4}
\begin{array}{rll}
 L_\xi(P_{i-1}\wedge P_i)
 \ggll
 H(P_{i})|P_{i-1}(\xi)|
 \ggll
 H(P_i)^{-\gamma},
 \quad \\ \\
 L_\xi(P_{i}\wedge P_{i+1})
 \ggll
 H(P_{i+1})|P_{i}(\xi)|
 \ggll
 H(P_{i+1})^{-\gamma}.
\end{array}
\end{equation}
Adjusting the implied constants if necessary, we may assume that
these estimates hold for all $i\in I$.

Consider now the next element $j=i_{k+1}$ of $I$.  By
construction, we have $V_i=V_{j-1}$ and so $P_i \wedge P_{i+1} =
\pm P_{j-1} \wedge P_j$.  Using (\ref{proof7-4}), we find
\[
 H(P_j)^{-\gamma}
 \ggll L_\xi(P_{j-1}\wedge P_j)
 = L_\xi(P_{i}\wedge P_{i+1})
 \ggll H(P_{i+1})^{-\gamma}
\]
and so, by (\ref{proof7-2}),
\begin{equation}
 \label{proof7-5}
 H(Q_{k+1})=H(P_j) \ggll H(P_{i+1}) \ggll H(Q_k)^\gamma.
\end{equation}
Moreover, since $P_i$ and $P_j$ are linearly independent (being
primary of distinct height), they form a basis of $V_i$.  Since
$V_i$ contains $P_i$ and $P_{i+1}$ whose resultant is non-zero, we
deduce that $P_i$ and $P_j$ also have a non-zero resultant and
using Lemma 2 of \cite{AR} together with (\ref{proof7-3}) and
(\ref{proof7-5}), we get
\[
 1
 \le |\Res(Q_k,Q_{k+1})|
 \ll H(Q_k) H(Q_{k+1})^2 |Q_k(\xi)|
 \ll 1.
\]
Finally, putting $\ell=i_{k+2}$, we observe that the polynomials
$P_i$, $P_j$ and $P_\ell$ are linearly independent as they span
the same vector space over $\bQ$ as $P_{j-1}$, $P_j$ and
$P_{j+1}$. Using Lemma 4 of \cite{AR} together with
(\ref{proof7-3}) and (\ref{proof7-5}), this gives
\[
 1
 \le |\det(Q_k,Q_{k+1},Q_{k+2})|
 \ll H(Q_{k+1}) H(Q_{k+2}) |Q_k(\xi)|
 \ll 1.
\]
The last three estimates together with (\ref{proof7-3}) show that
the sequence $(Q_k)_{k\ge 1}$ satisfies the condition (c).

To prove the last implication that (c) implies (b), assume the
existence of a sequence of polynomials $(Q_k)_{k\ge 1}$ as in
condition (c).  We first claim that $\xi$ is neither rational nor
quadratic over $\bQ$.  To prove this, assume on the contrary that
there exists a non-zero polynomial $Q\in\bZ[T]_{\le 2}$ which
vanishes at $\xi$.  Using Lemma 4 of \cite{AR}, we then find
\[
 |\det(Q,Q_k,Q_{k+1})|
 \ll H(Q)H(Q_{k+1})|Q_k(\xi)|
 \ll H(Q) H(Q_k)^{-\gamma^2},
\]
which implies that the integer $\det(Q,Q_k,Q_{k+1})$ is zero for
all sufficiently large values of $k$, against the hypothesis that
any three consecutive $Q_k$'s are linearly independent.  Define
\[
 \uy_k = Q_k \wedge Q_{k+1}
\]
for each index $k\ge 1$.  We have
\[
 \|\uy_k\| \ll H(Q_k)H(Q_{k+1}) \ll H(Q_k)^{\gamma^2}
\]
for all $k\ge 1$ and, applying Lemma \ref{lemma7-1} we deduce
\[
 L_\xi(\uy_k)
 \ggll H(Q_{k+1}) |Q_k(\xi)|
 \ggll H(Q_{k+1})^{-\gamma}
 \ll \|\uy_{k+1}\|^{-1/\gamma}.
\]
This shows in particular that $L_\xi(\uy_k)$ tends to zero as $k$
tends to infinity and thus that the sequence $(\uy_k)_{k\ge 1}$ is
unbounded.  So, for any sufficiently large real number $X$, there
exists an index $k$ such that $\|\uy_k\| \le X \le \|\uy_{k+1}\|$
and, for such a choice of $k$, the previous estimate gives
$L_\xi(\uy_k) \ll X^{-1/\gamma}$.  This shows that $\xi$ is
extremal.
\end{proof}

\appendix
\section{Finding new relations}

Let the notation be as in Section~\ref{apptriples}. Fix an integer
$k\ge 1$, a $(k+1)$-tuple of non-negative integers
$\ud=(d_0,\dots,d_k)$, an integer $p\ge 0$ and triples of
indeterminates $\uu_j = (u_{j,0},u_{j,1},u_{j,2})$ for
$j=0,1,\dots,k$.  To each monomial
\[
 \uu_0^{\ue_0}\cdots\uu_k^{\ue_k}
 = \prod_{j=0}^k\prod_{\ell=0}^2 u_{j,\ell}^{e_{j,\ell}}
 \in
 \bQ[\uu_0,\dots,\uu_k]
\]
we associate a {\it weight} given by
\[
 \sum_{j=0}^k\sum_{\ell=0}^2 \ell e_{j,\ell}
\]
and a {\it multi-degree} given by
\[
 (|\ue_0|,\dots,|\ue_k|)
  =(e_{0,0}+e_{0,1}+e_{0,2},\dots,e_{k,0}+e_{k,1}+e_{k,2}).
\]
We denote by $E(\ud,p)$ the sub-$\bQ$-vector space of
$\bQ[\uu_0,\dots,\uu_k]$ generated by all monomials of
multi-degree $\ud$ and weight $p$. We also denote by $\uun$ the
element of $\bQ^{3k+3}$ all of whose coordinates are equal to $1$.

\begin{proposition}
  \label{searchpol}
  Suppose that a polynomial $P\in E(\ud,p)$ satisfies
\[
 \left(
 \prod_{j=0}^k\prod_{\ell=1}^2
 \left(\frac{\partial}{\partial u_{j,\ell}}\right)^{f_{j,\ell}}
 P
 \right)(\uun)
 = 0
\]
for any choice of non-negative integers $f_{0,1}, f_{0,2}, \dots,
f_{k,1}, f_{k,2}$ with
\[
2\sum_{j=0}^k (f_{j,1}+f_{j,2})\gamma^j
  \le \sum_{j=0}^k d_j\gamma^j
  \et
f_{j,1}+f_{j,2} \le d_j,
  \quad (0\le j\le k).
\]
Then, we have
\[
P(\uy_i,\uy_{i+1},\dots,\uy_{i+k}) = 0
\]
for all sufficiently large values of $i$.
\end{proposition}

\begin{proof}
We find
\[
\begin{array}{rl}
P(\uy_i,\dots,\uy_{i+k})
 &=
 \disp
 \xi^p y_{i,0}^{d_0}\cdots y_{i+k,0}^{d_k}
 P\left(
    1,
    \frac{y_{i,1}}{y_{i,0}\xi},
    \frac{y_{i,2}}{y_{i,0}\xi^2},
    \dots,
    1,
    \frac{y_{i+k,1}}{y_{i+k,0}\xi},
    \frac{y_{i+k,2}}{y_{i+k,0}\xi^2}
    \right) \\
 &=
 \disp
 \xi^p y_{i,0}^{d_0}\cdots y_{i+k,0}^{d_k}
 \sum
 \prod_{j=0}^k \prod_{\ell=1}^2
 \frac{1}{f_{j,\ell}!}
 \left(
  \left(
  \frac{y_{i+j,\ell}}{y_{i+j,0}\xi^\ell}-1
  \right)
  \frac{\partial}{\partial u_{j,\ell}}
 \right)^{f_{j,\ell}}
 P (\uun), \\
\end{array}
\]
where the sum in the second expression ranges over all choices of
non-negative integers $f_{0,1},f_{0,2},\dots,f_{k,0},f_{k,1}$ with
$f_{j,1}+f_{j,2}\le d_j$ for $j=0,\dots,k$.  Since, for any such
choice of integers, we have
\[
\begin{array}{rl}
 \disp
 \left|
 y_{i,0}^{d_0}\cdots y_{i+k,0}^{d_k}
 \prod_{j=0}^k \prod_{\ell=1}^2
  \left(
  \frac{y_{i+j,\ell}}{y_{i+j,0}\xi^\ell}-1
  \right)^{f_{j,\ell}}
 \right|
 &\disp
 \ll
 \prod_{j=0}^k \|\uy_{i+j}\|^{d_j-2f_{j,1}-2f_{j,2}} \\
 &\disp
 \ll
 \|\uy_i\|^{\sum_{j=0}^k (d_j-2f_{j,1}-2f_{j,2})\gamma^j}, \\
\end{array}
\]
the hypothesis implies
\[
 |P(\uy_i,\dots,\uy_{i+k})|
 \ll
 \|\uy_i\|^{-\epsilon}
\]
for some positive real number $\epsilon$.  Since $P(\uy_i,\dots,
\uy_{i+k})$ is a rational number with bounded denominator, it must
therefore vanish for all sufficiently large values of $i$.
\end{proof}

Empirically it seems that, for a given multi-degree $\ud$, the
dimension of $E(\ud,p)$ is maximal with $p=|\ud|=d_0+\dots+d_k$.
For values of $p$ at equal distance from $|\ud|$, that is for
integers $p,q\ge 0$ with $p+q = 2|\ud|$, the dimensions of the
corresponding vector spaces $E(\ud,p)$ and $E(\ud,q)$ are the
same.  So, it is natural to look first for polynomials in
$E(\ud,|\ud|)$.

\medskip
A computer search based on the above proposition found two
non-zero polynomials with the appropriate vanishing, namely
\[
 \det(\uu_0,\uu_1,[\uu_3,\uu_3,\uu_4])
 \et
 \det(\uu_1,\uu_2,[\uu_3,\uu_3,\uu_4])
\]
where the symbol $[\uu_3,\uu_3,\uu_4]$ is defined by
(\ref{bracket}).  The first polynomial has multi-degree
$(1,1,0,2,1)$ and weight $5$, while the other has multi-degree
$(0,1,1,2,1)$ and weight $5$. They provide the relations
(\ref{relations}) of Proposition~\ref{prop-rel}.


\end{document}